# Solution of Fuzzy Growth and Decay Model

U. M. Pirzada

School of Engineering and Technology,
Navrachana University of Vadodara,
salmap@nuv.ac.in

**Abstract:** Mathematical modelling for population growth leads to a differential equation. In population growth model, we assume that rate increase of population is proportional to current population. That is, $dx/dt = kx$, x is a current population, k is proportionality constant represents growth rate. But in real situation, it is often ambiguous to determine exact amount of current population. It can be measured approximately. For instance, initially number of bacteria is approximately 20 and this approximate number can be represented using fuzzy number. Therefore, the appropriate growth or decay model is described using fuzzy concept. With this motivation, this paper presents solution of fuzzy growth and decay model. The solution is analysed using Seikkala differentiability of fuzzy-valued function.

**1 Preliminary**

**Definition 1.1** Let R be the set of real numbers and $a : R \rightarrow [0, 1]$ be a fuzzy set. We say that a is a fuzzy number if it satisfies the following properties:

(i) a is normal, that is, there exists $r_0 \in R$ such that $a(r_0) = 1$;

(ii) a is fuzzy convex, that is, $a(tr + (1 − t)s) \geq \min\{a(r), a(s)\}$, whenever $r, s \in R$ and $t \in [0, 1]$;

(iii) $a(r)$ is upper semi-continuous on R, that is, $\{r / a(r) \geq \alpha\}$ is a closed subset of R for each $\alpha \in (0, 1]$;

(iv) $cl\{r \in R / a(r) > 0\}$ forms a compact set,

where cl denotes closure of a set. The set of all fuzzy numbers on R is denoted by F(R).

**Definition 1.2** For all $\alpha \in (0, 1]$, $\alpha$-level set $a_\alpha$ of any $a \in F(R)$ is defined as $a_\alpha = \{r \in R \,/\, a(r) \geq \alpha\}$.

The 0-level set $a_0$ is defined as the closure of the set $\{r \in R \,/\, a(r) > 0\}$.

**Remarks:**
  (i) By definition of fuzzy numbers, it can be prove that, for any $a \in F(R)$ and for each $\alpha \in (0, 1]$, $a_\alpha$ is compact convex subset of R, and we write $a_\alpha = [a_1(\alpha), a_2(\alpha)]$.
  (ii) The fuzzy number $a \in F(R)$ can be recovered from its $\alpha$-level sets by a well-known decomposition theorem (ref. [1]), which states that $a = \bigcup_{\alpha \in [0,1]} (\alpha a_\alpha)$, where union on the right-hand side is the standard fuzzy union.

The following Theorem of Goetschel and Voxman [2], shows the characterization of a fuzzy number in terms of its $\alpha$-level sets.

**Theorem 1.1** For $a \in F(R)$, define two functions $a_1(\alpha), a_2(\alpha) : [0, 1] \in R$. Then
(i) $a_1(\alpha)$ is bounded left continuous non-decreasing function on $(0,1]$;
(ii) $a_2(\alpha)$ is bounded left continuous non-increasing function on $(0,1]$;
(iii) $a_1(\alpha)$ and $a_2(\alpha)$ are right continuous at $\alpha = 0$;
(iv) $a_1(\alpha) \leq a_2(\alpha)$.

Moreover, if the pair of functions $a_1(\alpha)$ and $a_2(\alpha)$ satisfy the conditions (i)-(iv), for each $\alpha \in [0, 1]$, then there exists a unique $a \in F(R)$ such that $a_\alpha = [a_1(\alpha), a_2(\alpha)]$, for $\alpha \in [0, 1]$.

**Definition 1.3** The membership function of a triangular fuzzy number a is defined as
$$a(r) = \begin{cases} \frac{r-a^L}{a-a^L} & if\, a^L \leq r \leq a \\ \frac{a^U-r}{a^U-a} & if\, a < r \leq a^U \\ 0 & otherwise \end{cases}$$

and it is denoted by $a = (a^L, a, a^U)$. The $\alpha$-level sets of a are

$a_\alpha = [(1 - \alpha)a^L + \alpha a, (1 - \alpha)a^U + \alpha a]$.

**Definition 1.4** According to Zadeh's extension principle, addition, multiplication of two fuzzy numbers a, b and scalar multiplication of fuzzy number with a scalar $\lambda \in R$ by their α-level sets are defined as follows:

$(a + b)_\alpha = [a_1(\alpha) + b_1(\alpha), a_2(\alpha) + b_2(\alpha)]$,

$(a \times b)_\alpha = [\min\{a_1(\alpha) b_1(\alpha), a_1(\alpha) b_2(\alpha), a_2(\alpha) b_1(\alpha), a_2(\alpha) b_2(\alpha)\}, \max\{a_1(\alpha) b_1(\alpha), a_1(\alpha) b_2(\alpha), a_2(\alpha) b_1(\alpha), a_2(\alpha) b_2(\alpha)]$,

$(\lambda a)_\alpha = [\lambda a_1(\alpha), \lambda a_2(\alpha)]$, if $\lambda \geq 0$
$\phantom{(\lambda a)_\alpha} = [\lambda a_2(\alpha), \lambda a_1(\alpha)]$, if $\lambda < 0$,

where α-level sets of a and b are $a_\alpha = [a_1(\alpha), a_2(\alpha)]$, $b_\alpha = [b_1(\alpha), b_2(\alpha)]$, for $\alpha \in [0, 1]$.

**Definition 1.5** Let V be a real vector space and F(R) be a set of fuzzy numbers. Then a fuzzy-valued function $y : V \to F(R)$ is defined on V. Corresponding to such a function y and for each $\alpha \in [0, 1]$, we denote two real-valued functions $y_1(t, \alpha)$ and $y_2(t, \alpha)$ on V for all $t \in V$. These functions $y_1(t, \alpha)$ and $y_2(t, \alpha)$ are called α-level functions of the fuzzy-valued function y.

Seikkala differentiability of fuzzy-valued function y is defined as follows. We adopted the definition from [3].

**Definition 1.6** Let I be subsets of R. Let y be a fuzzy-valued function defined on I. Let α-level sets $y_\alpha(t) = [y_1(t,\alpha), y_2(t,\alpha)]$ for all α. We assume that $y_i(t,\alpha)$ have derivatives, for all $t \in I$, for each α, i = 1, 2.

Define $(y'(t))_\alpha = [y_1'(t, \alpha), y_2'(t, \alpha)]$ for all $t \in I$, all α.

If, for each fixed $t \in I$, $(y'(t))_\alpha$ defines the α-level set of a fuzzy number, then we say that y(t) is differentiable at t.

The sufficient conditions for $(y'(t))_\alpha$ to define α-level sets of a fuzzy numbers are

(i) $y_1'(t, \alpha)$ is an increasing function of α for each $t \in I$;
(ii) $y_2'(t, \alpha)$ is a decreasing function of α for each $t \in I$; and
(iii) $y_1'(t, \alpha) \leq y_2'(t, \alpha)$ for all $t \in I$.

**2 Fuzzy growth and decay model**

Derivation of model for population growth leads to a differential equation. In population growth model, we assume that rate increase of population is proportional to current population. That is, dx / dt = kx, x is a current population, k is proportionality constant represents growth rate. For population of bacteria, we cannot say at certain time amount of bacteria is exactly 10 or 20, it can be measured approximately. That is, initially number of bacteria is approximately 20. So this approximate number can be represented using fuzzy number. In solution of growth model, we cannot get exact amount of population. The appropriate growth model is described using fuzzy concept.

Fuzzy growth and decay models are studied by Buckley and Feuring in [4] as an application of fuzzy differential equations. S. P. Mondal et al.[5] have also studied solution of fuzzy growth and decay model. They have not considered any differentiability concept to study the model. In this paper, solution of fuzzy growth and decay model is studied using Seikkala differentiability of fuzzy-valued function.

Consider a fuzzy initial value problem (FIVP)

$dy/dt = f(t, y) = k \times y$, $y(0) = c$,

where $y: I \to F(R)$ is a unknown fuzzy-valued function, for $t \in I$, k and c are fuzzy numbers and $f: I \times F(R) \to F(R)$ is a fuzzy-valued function. This (FIVP) is a fuzzy growth / decay model.

We assume that fuzzy-valued function y is Seikkala differentiable on I. Using fuzzy arithmetic, the (FIVP) can be written as system of crisp parametric differential equations

$dy_1/dt = f_1(t, y_1, y_2, \alpha) = \min\{k_1(\alpha) y_1(t, \alpha), k_1(\alpha) y_2(t, \alpha), k_2(\alpha) y_1(t, \alpha), k_2(\alpha) y_2(t, \alpha)\}$ **(1)**

$dy_2/dt = f_2(t, y_1, y_2, \alpha) = \max\{k_1(\alpha) y_1(t, \alpha), k_1(\alpha) y_2(t, \alpha), k_2(\alpha) y_1(t, \alpha), k_2(\alpha) y_2(t, \alpha)\}$ **(2)**

$y_1(0, \alpha) = c_1(\alpha)$ and $y_2(0, \alpha) = c_2(\alpha)$, for $t \in I$ and $\alpha$.

## 2.1 Solution

To find solution of fuzzy growth and decay model, we consider two cases.

**Case I:** $k_1(\alpha), k_2(\alpha) \geq 0$. Fuzzy growth model

We also assume $y_1(t, \alpha)$ and $y_2(t, \alpha)$ are positive. Therefore, the above system can be simplified as

$dy_1/dt = f_1(t, y_1, y_2, \alpha) = k_1(\alpha) y_1(t, \alpha)$,

$dy_2/dt = f_2(t, y_1, y_2, \alpha) = k_2(\alpha) y_2(t, \alpha)$, $y_1(0, \alpha) = c_1(\alpha)$ and $y_2(0, \alpha) = c_2(\alpha)$, for $t \in I$ and $\alpha$.

Solving this system, we get

$y_1(t, \alpha) = c_1(\alpha) \exp(k_1(\alpha)t)$, $y_2(t, \alpha) = c_2(\alpha) \exp(k_2(\alpha)t)$ for $t \in I$ and $\alpha$.

We assume that $c_1(\alpha)$, $c_2(\alpha) \geq 0$. Now we will check $[y_1(t, \alpha), y_2(t, \alpha)]$ defines fuzzy number or not for each t in I. Moreover, we need to check $[y_1'(t, \alpha), y_2'(t, \alpha)]$ defines fuzzy number for each t in I so that we say that solution y is Seikkala differentiable.

That is, we need to check that $\partial y_1(t, \alpha) / \partial \alpha > 0$ and $\partial y_2(t, \alpha) / \partial \alpha < 0$. We see that

$\partial y_1(t, \alpha) / \partial \alpha = c_1'(\alpha) \exp(k_1(\alpha)t) + c_1(\alpha) k_1(\alpha) \exp(k_1(\alpha)t) k_1'(\alpha) > 0$

as $c_1'(\alpha) > 0$ and $k_1'(\alpha) > 0$

$\partial y_2(t, \alpha) / \partial \alpha = c_2'(\alpha) \exp(k_2(\alpha)t) + c_2(\alpha) k_2(\alpha) \exp(k_2(\alpha)t) k_2'(\alpha) < 0$

as $c_2'(\alpha) < 0$ and $k_2'(\alpha) < 0$.

Hence we say that $[y_1(t, \alpha), y_2(t, \alpha)]$ defines fuzzy number for each t in I. We also need to check that $[y_1'(t, \alpha), y_2'(t, \alpha)]$ defines fuzzy number for each t in I.

Since $y_1'(t, \alpha) = c_1(\alpha) k_1(\alpha) \exp(k_1(\alpha)t)$, $y_2'(t, \alpha) = c_2(\alpha) k_1(\alpha) \exp(k_2(\alpha)t)$, it can easily seen that

$\partial y_1'(t, \alpha) / \partial \alpha > 0$ and $\partial y_2'(t, \alpha) / \partial \alpha < 0$ and therefore, y is Seikkala differentiable on I. Hence, $y = c \exp(kt)$, where c and k are fuzzy numbers, is a fuzzy solution of fuzzy growth model on domain I.

**Case II:** $k_1(\alpha), k_2(\alpha) < 0$. Fuzzy decay model

We assume $y_1(t, \alpha)$ and $y_2(t, \alpha)$ are positive. Equations (1) and (2) can written as

$dy_1 / dt = f_1(t, y_1, y_2, \alpha) = k_1(\alpha) y_2(t, \alpha)$,

$dy_2 / dt = f_2(t, y_1, y_2, \alpha) = k_2(\alpha) y_1(t, \alpha)$,

with $y_1(0, \alpha) = c_1(\alpha)$ and $y_2(0, \alpha) = c_2(\alpha)$, for $t \in I$ and $\alpha$.

We assume that $c_1(\alpha), c_2(\alpha) \geq 0$. The solution of above system of equations (Adopted from Buckley's paper [4])

$y_1(t, \alpha) = A_{11}(\alpha) \exp(pt) + A_{12}(\alpha) \exp(-pt)$; $y_2(t, \alpha) = A_{21}(\alpha) \exp(pt) - A_{22}(\alpha) \exp(-pt)$,

where $p = \sqrt{k_1(\alpha)k_2(\alpha)}$, $q = \sqrt{k_1(\alpha)/k_2(\alpha)}$, $A_{11}(\alpha) = 0.5(c_1(\alpha) + q\, c_2(\alpha))$; $A_{12}(\alpha) = 0.5(c_1(\alpha) - q\, c_2(\alpha))$; $A_{21}(\alpha) = 0.5(c_1(\alpha) / q + c_2(\alpha))$ and $A_{22}(\alpha) = 0.5(c_1(\alpha) / q - c_2(\alpha))$ for $t \in I$ and $\alpha$.

We need to check [$y_1(t, \alpha)$, $y_2(t, \alpha)$] defines fuzzy number or not for each t in I. Moreover, we need to check [$y_1'(t, \alpha)$, $y_2'(t, \alpha)$] defines fuzzy number for each t in I so that we say that solution y is Seikkala differentiable.

For this we need to check sufficient conditions for existence of fuzzy numbers. That is, we need to check that $\partial y_1(t, \alpha) / \partial \alpha > 0$ and $\partial y_2(t, \alpha) / \partial \alpha < 0$.

We also need to check that check [$y_1'(t, \alpha)$, $y_2'(t, \alpha)$] defines fuzzy number for each t in I. That is, $\partial y_1'(t, \alpha) / \partial \alpha > 0$ and $\partial y_2'(t, \alpha) / \partial \alpha < 0$ so that y is Seikkala differentiable on I. Then y is a fuzzy solution of fuzzy decay model.

To analyse the model, we consider specific values for k and c. Let k = -1 and c = (2, 4, 6) be a triangular fuzzy number. Then system can be re-written as

$dy_1 / dt = f_1(t, y_1, y_2, \alpha) = -y_2(t, \alpha)$,

$dy_2 / dt = f_2(t, y_1, y_2, \alpha) = -y_1(t, \alpha)$,

with $y_1(0, \alpha) = (2 + 2\alpha)$ and $y_2(0, \alpha) = (6 - 2\alpha)$, for $t \in I$ and $\alpha$. Here $k_1(\alpha) = k_2(\alpha) = -1$.

The solution of the system is

$y_1(t, \alpha) = A_{11}(\alpha) \exp(pt) + A_{12}(\alpha) \exp(-pt)$; $y_2(t, \alpha) = A_{21}(\alpha) \exp(pt) - A_{22}(\alpha) \exp(-pt)$,

where $p = \sqrt{k_1(\alpha) k_2(\alpha)} = 1$, $q = \sqrt{k_1(\alpha) / k_2(\alpha)} = 1$, $A_{11}(\alpha) = 0.5(c_1(\alpha) + c_2(\alpha))$; $A_{12}(\alpha) = 0.5(c_1(\alpha) - c_2(\alpha))$; $A_{21}(\alpha) = 0.5(c_1(\alpha) + c_2(\alpha))$ and $A_{22}(\alpha) = 0.5(c_1(\alpha) - c_2(\alpha))$ for $t \in I$ and $\alpha$.

We need to check [$y_1(t, \alpha)$, $y_2(t, \alpha)$] defines fuzzy number or not for each t in I. That is, we need to check that $\partial y_1(t, \alpha) / \partial \alpha > 0$ and $\partial y_2(t, \alpha) / \partial \alpha < 0$.

$\partial y_1(t, \alpha) / \partial \alpha = A_{11}'(\alpha) \exp(pt) + A_{11}(\alpha) t \exp(pt) p'(\alpha) + A_{12}'(\alpha) \exp(-pt) + A_{12}(\alpha) (-t) \exp(-pt) p'(\alpha)$.

But $p'(\alpha) = 0$ as $p(\alpha) = 1$, we have

$\partial y_1(t, \alpha) / \partial \alpha = A_{11}'(\alpha) \exp(pt) + A_{12}'(\alpha) \exp(-pt)$

$= 0.5(c_1'(\alpha) + c_2'(\alpha)) \exp(pt) + 0.5(c_1'(\alpha) - c_2'(\alpha)) \exp(-pt)$.

Using $c_1'(\alpha) = 2$ and $c_2'(\alpha) = -2$,

$= 0.5(2 - 2) \exp(pt) + 0.5(2 - (-2)) \exp(-pt)$

$= 2 \exp(-t) = 2 / \exp(t) > 0$ for $t > 0$ in I.

Now,

$\partial y_2(t, \alpha) / \partial \alpha = A_{21}'(\alpha) \exp(pt) - A_{22}'(\alpha) \exp(-pt)$

$\quad = 0.5(c_1'(\alpha) + c_2'(\alpha)) \exp(pt) - 0.5(c_1'(\alpha) - c_2'(\alpha)) \exp(-pt)$

$\quad = 0.5(2 - 2) \exp(pt) - 0.5(2 - (-2)) \exp(-pt)$

$\quad = -2 \exp(-t) = 2/ \exp(t) < 0$ for $t > 0$ in I.

Therefore $[y_1(t, \alpha), y_2(t, \alpha)]$ defined fuzzy number for each t in I.

Further, we need to check $[y_1'(t, \alpha), y_2'(t, \alpha)]$ defines fuzzy number for each t in I so that we say that solution y is Seikkala differentiable. That is, we need to check

$\partial y_1'(t, \alpha) / \partial \alpha > 0$ and $\partial y_2'(t, \alpha) / \partial \alpha < 0$.

Consider,

$y_1'(t, \alpha) = A_{11}(\alpha) \exp(pt) p - A_{12}(\alpha) \exp(-pt) p$ and

$y_2'(t, \alpha) = A_{21}(\alpha) \exp(pt) p + A_{22}(\alpha) \exp(-pt) p$.

Now,

$\partial y_1'(t, \alpha) / \partial \alpha = A_{11}'(\alpha) \exp(pt) - A_{12}'(\alpha) \exp(-pt)$

$\quad = 0.5(c_1'(\alpha) + c_2'(\alpha)) \exp(pt) - 0.5(c_1'(\alpha) - c_2'(\alpha)) \exp(-pt)$, but $c_1'(\alpha) = 2$ and $c_2'(\alpha) = -2$,

$\quad = 0.5(2 - 2) \exp(pt) - 0.5(2 - (-2)) \exp(-pt)$

$\quad = -2 \exp(-t) = 2/ \exp(t) < 0$ for $t > 0$ in I

Similar way, we see that $\partial y_2'(t, \alpha) / \partial \alpha > 0$ for $t > 0$ in I.

Therefore, $[y_1'(t, \alpha), y_2'(t, \alpha)]$ does not define fuzzy number for t in I and hence we say that solution y is not Seikkala differentiable. Therefore, we conclude that fuzzy decay model has no solution under Seikkala differentiability of fuzzy-valued function.

# 3 Conclusions

Solution of fuzzy growth and decay model is studied using Seikkala differentiability of fuzzy-valued function. The solution is analysed for different values of k. We observe that solution of fuzzy decay model with crisp coefficient k = -1, does not exist under Seikkala differentiability, while solution of fuzzy growth model exists with fuzzy coefficient k and fuzzy initial conditions. In future, we would like to study the solution of fuzzy decay model under generalized differentiability concept.